\documentclass{ifacconf}

\usepackage{graphicx}      
\usepackage{natbib}        

\usepackage{amssymb} 
\usepackage{amsmath}
\usepackage{color}
\usepackage{graphicx,pst-all}

\newcommand{\cD}{\mathcal{D}}

\newcommand{\N}{{\mathbb{N}}}

\newcommand{\R}{{\mathbb{R}}}

\DeclareMathOperator{\im}{im}

\DeclareMathOperator{\rf}{ref}

\begin{document}
\begin{frontmatter}

\title{Learning-based Funnel-MPC for output-constrained nonlinear systems\thanksref{footnoteinfo}}

\thanks[footnoteinfo]{Thomas Berger and Karl Worthmann gratefully acknowledge support from the German Research Foundation (Deutsche Forschungsgemeinschaft) via the grants BE\,6263/1-1 and WO\,2056/6-1.}

\author[First]{Thomas Berger}
\author[Second]{Carolin K\"{a}stner}
\author[Second]{Karl Worthmann}

\address[First]{Institut f\"ur Mathematik, Universit\"at Paderborn, Warburger Str.~100, 33098~Paderborn, Germany (e-mail: thomas.berger@math.upb.de).}
\address[Second]{Institute for Mathematics, 
	Technische Universit\"{a}t Ilmenau, Ilmenau, Germany (e-mail: \{carolin.kaestner, karl.worthmann\}@tu-ilmenau.de)}

\begin{abstract}                
	We exploit an adaptive control technique, namely funnel control, %
	in order to establish both initial and recursive feasibility in %
	Model Predictive Control (MPC) for output-constrained nonlinear systems. 
	Moreover, we show that the resulting feedback controller outperforms %
	the funnel controller both w.r.t.\ the required sampling rate %
	for a zero-order-hold implementation and required control action. %
	We further propose a combination of funnel control and MPC, %
	exploiting the performance guarantees of the model-free funnel controller %
	during a learning phase and the advantages of the model-based MPC scheme thereafter.
\end{abstract}

\begin{keyword}
	model predictive control, funnel control,  learning, output tracking, initial feasibility, parameter identification, recursive feasibility
\end{keyword}

\end{frontmatter}

\section{Introduction}

Model Predictive Control (MPC) is nowadays a widely-applied control technique, %
largely thanks to its applicability to constrained nonlinear %
multi-input, multi-output systems, see e.g.~\cite{GruePann17} or~\cite{RawlMayn18}. %
However, there are two main obstacles: %
On the one hand, MPC requires a model for its prediction and optimization step, and on the other %
initial and recursive feasibility have to be ensured. %
In the present paper, we propose learning-based Funnel-MPC in order to resolve these issues. %
To this end, we exploit the concept of funnel control --~a model-free output-error feedback of high-gain type developed in~\cite{IlchRyan02b}, see also the survey~\cite{IlchRyan08}.\\
MPC requires either a sufficiently long prediction horizon (see, e.g.~\cite{BoccGrun14}) or %
suitably constructed terminal constraints (see, e.g.~\cite{RawlMayn18}) in order to guarantee recursive feasibility while initial feasibility is assumed. %
Both approaches and their respective prerequisites 
are difficult to achieve in the presence of (time-varying) state (or output) constraints.\\
The complementary concept of funnel control guarantees a prescribed tracking performance over the whole time interval.
The funnel controller proved its potential for tracking problems in various applications, see e.g.~\cite{BergRaue18,BergReis14a,Hack17} and the references therein. We investigate the combination of MPC and funnel control in order to benefit from the best of both worlds: %
guaranteed feasibility (funnel control) and superior performance (MPC).\\
The key idea is reflected by the ``funnel-like'' cost function, %
which is based on the model-free control law used in funnel control and becomes infinite when the tracking error approaches the funnel boundary. %
An immediate advantage of this penalization-based approach is %
that available results in funnel control guarantee the existence of a control input which meets the output constraints. %
We exploit this to prove initial and recursive feasibility of a corresponding MPC scheme.\\
However, MPC still requires a model of the system while funnel control achieves the objective with structural information only, without any knowledge of specific system parameters. To resolve this, we additionally employ learning techniques to obtain the model parameters. Starting with a raw estimate, the initial state and the system dynamics are approximated during a learning phase while adherence of the output constraints is guaranteed by the funnel controller. %
After that, the mechanism switches to the Funnel-MPC scheme, thus improving controller performance.\\
We like to emphasize that the Funnel-MPC scheme also significantly relaxes requirements on the sampling rate, while still considerably improving the performance and the range of applied control values compared to (a zero-order-hold implementation of) the funnel controller.

\noindent\textbf{Notation}: $\mathbb{N}$ and $\mathbb{R}$ denote the natural and real numbers, resp., $\N_0 = \mathbb{N} \cup \{0\}$ and $\mathbb{R}_{\geq 0} = [0,\infty)$. $\mathcal{C}^k(\mathbb{R}_{\geq 0}, \mathbb{R})$ is the (linear) space of $k$-times continuously differentiable functions $f:\R_{\ge 0}\to\R$, while $\mathcal{L}^\infty_{\operatorname{loc}}(\mathbb{R}_{\geq 0}, \mathbb{R})$ denotes the (linear) space of Lebesgue-measurable, and locally essentially bounded functions. Moreover, we define $\mathcal{B}^{k,\infty}(\mathbb{R}_{\geq 0},\mathbb{R})$ as the space of $k$-times continuously differentiable functions with $f,\dot f,\ldots, f^{(k)}\in \mathcal{L}^\infty(\mathbb{R}_{\geq 0}, \mathbb{R})$.

\section{Combining Funnel Control and MPC}

We consider the control affine system
\begin{equation}\label{eq:Sys}
\begin{aligned}
	\dot{x}(t) & = f(x(t)) + g(x(t)) u(t),\\
	y(t) & = h(x(t))
\end{aligned}
\end{equation}
with sufficiently smooth vector fields $f: \mathbb{R}^n \rightarrow \mathbb{R}^n$, $g: \mathbb{R}^{n} \rightarrow \mathbb{R}^{n \times p}$, a sufficiently smooth mapping $h: \mathbb{R}^n \rightarrow \mathbb{R}^p$, and control input function $u \in \mathcal{L}^\infty_{\operatorname{loc}}(\mathbb{R}_{\geq 0}, \mathbb{R}^p)$. Note that the dimensions of the output and input coincide.\\
For single-input, single-output (SISO) systems, i.e., $p = 1$, the control affine system \eqref{eq:Sys} is said to have relative degree $r \in \mathbb{N}$, if the conditions
\begin{align*}
	\forall\;k \in \{1,\ldots,r-1\}:\ L_g L_f^{k-1} h(x) & = 0, \\
	L_g L_f^{r-1} h(x) & \neq 0
\end{align*}
hold for all $x \in \mathbb{R}^n$, see e.g.~\cite{Isid95}. %
Recall that the \emph{Lie derivative} of~$h$ along~$f$ is defined by
\[
    \big(L_f h\big)(x) = \sum_{i=1}^n \frac{\partial h}{\partial x_i}(x)\, f_i(x) = h'(x) f(x),
\]
and we may successively define $L_f^k h = L_f (L_f^{k-1} h)$ with $L_f^0 h = h$.\\
For multi-input, multi-output (MIMO) systems, i.e., $p \geq 2$, again we have $(L_f h)(x) = h'(x) f(x)$, where $h'$ is the Jacobian of~$h$. Furthermore, denoting with $g_i(x)$ the columns of~$g(x)$ for $i=1,\ldots,p$, we have
\[
    \big(L_g h\big)(x) = [\big(L_{g_1} h\big)(x), \ldots, \big(L_{g_p} h\big)(x)].
\]
Then the system~\eqref{eq:Sys} has relative degree $r \in \mathbb{N}$, if
\[
	\forall\;k \in \{1,\ldots,r-1\}\ \forall\, x\in\R^n:\ L_g L_f^{k-1} h(x)  = 0
\]
and the high-gain matrix
\begin{equation}
	\Gamma(x) := \big(L_g L_f^{r-1} h\big)(x) \label{eq:Gamma}
\end{equation}
is invertible for all~$x\in\R^n$. Note that, for simplicity, we require the above properties to hold for all $x\in\R^n$, while usually these are \emph{local} properties, i.e., the high-gain matrix may only be invertible on different separated open subsets of~$\R^n$.\\
If~\eqref{eq:Sys} has relative degree $r\in\N$, then there exists a diffeomorphism $\Phi:\R^n\to\R^n$ such that the coordinate transformation $(\xi(t),\eta(t)) = \Phi(x(t))$, where $\xi(t) = (y(t), \dot y(t),\ldots,y^{(r-1)}(t))$, puts the system into Byrnes-Isidori form (cf.~\cite[Sec.~5.1]{Isid95})
\begin{equation}\label{eq:BIF}
\begin{aligned}
    y^{(r)}(t) &= p\big(\xi(t),\eta(t)\big) + \Gamma\big(\Phi^{-1}\big(\xi(t),\eta(t)\big)\big)\, u(t),\\
    \dot \eta(t) &= q\big(\xi(t),\eta(t)\big).
\end{aligned}
\end{equation}
If $p\ge 2$, then we need to require that the distribution $G(x) = \im g(x)$ is involutive in order for this to be feasible. We call $G(x)$ involutive (see e.g.~\cite[Sec.~1.3]{Isid95}), if for all smooth vector fields $g_1, g_2: \R^n \to \R^n$ with $g_i(x)\in G(x)$ for all $x \in \R^n$ and $i=1,2$ we have that the Lie bracket $[g_1,g_2](x) = g_1'(x) g_2(x) - g_2'(x) g_1(x)$ satisfies $[g_1,g_2](x) \in G(x)$ for all $x\in\R^n$.

\subsection{Funnel control revisited}

\begin{figure}[htb]
\hspace{6mm}\includegraphics[width=8cm, trim=170 550 230 120, clip]{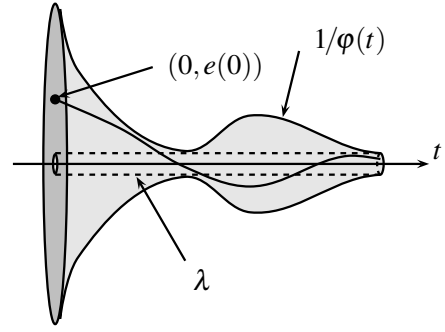}
	\caption{Error evolution in a funnel with boundary~$1/\varphi(t)$.}
	\label{Fig:funnel}
\end{figure}

For a given reference trajectory $y_{\operatorname{ref}} \in \mathcal{B}^{r,\infty}(\mathbb{R}_{\geq 0},\mathbb{R}^p)$ %
the objective in funnel control is to design a control law such that the tracking error~$e(t) = y(t) - y_{\rf}(t)$ evolves within a performance funnel determined by a function $\varphi \in \mathcal{B}^{r,\infty}(\mathbb{R}_{\geq 0},\mathbb{R})$ which satisfies $\varphi(t)>0$ for all $t>0$ and $\liminf_{t\rightarrow \infty} \varphi(t) > 0$,  see Fig.~\ref{Fig:funnel}. Furthermore, all signals $u, e, \dot e, \ldots, e^{(r-1)}$ in the closed-loop system should remain bounded.\\
The boundary of the performance funnel is given by~$1/\varphi$. If we choose $\varphi(0)=0$, which is explicitly allowed, then the initial error can be arbitrary since $\varphi(0) \|e(0)\| < 1$; in this case the funnel boundary $1/\varphi$ has a pole at $t=0$. Since $\varphi$ is bounded, there exists $\lambda>0$ such that $1/\varphi(t) \ge \lambda$ for all $t>0$.\\
It was shown in~\cite{BergLe18} that under some structural assumptions the funnel controller
\begin{equation}\label{eq:fun-con}
\boxed{\begin{aligned}
u_{\operatorname{FC}}(t) &= \sigma\, k_{r-1}(t) \, e_{r-1}(t),\\
e_0(t)&=e(t) = y(t) - y_{\rm ref}(t),\\
e_1(t)&=\dot{e}_0(t)+k_0(t)\,e_0(t),\\
e_2(t)&=\dot{e}_1(t)+k_1(t)\,e_1(t),\\
& \ \ \vdots \\
e_{r-1}(t)&=\dot{e}_{r-2}(t)+k_{r-2}(t)\,e_{r-2}(t),\\
k_i(t)&=\frac{1}{1-\varphi_i(t)^2\|e_i(t)\|^2},\quad i=0,\dots,r-1,
\end{aligned}
}
\end{equation}
achieves the above described control objective (more precisely, $\|e_i(t)\| < \varphi_i(t)^{-1}$ for all $t>0$), where $\sigma\in\{-1,1\}$ and each $\varphi_k \in \mathcal{B}^{r-k,\infty}(\mathbb{R}_{\geq 0},\mathbb{R})$ satisfies $\varphi_k(t)>0$ for all $t>0$ and $\liminf_{t \rightarrow \infty} \varphi_k(t) > 0$, $k=0,\ldots,r-1$.
The structural assumptions required in~\cite{BergLe18} are knowledge of the relative degree~$r$, a bounded-input, bounded-output (BIBO) property of the internal dynamics (i.e., the second equation in~\eqref{eq:BIF}), and sign-definite high-gain matrix, i.e.,
\[
    \forall\, x\in\R^n\ \forall\, v\in\R^p:\quad v^\top \Gamma(x) v = 0\ \ \iff\ \ v=0;
\]
this is equivalent to either $\Gamma(x)+\Gamma(x)^\top > 0$ (in which case we have $\sigma = -1$ in~\eqref{eq:fun-con}) or  $\Gamma(x)+\Gamma(x)^\top < 0$ (in which case we have $\sigma = 1$ in~\eqref{eq:fun-con}).\\
If the system~\eqref{eq:Sys} and the current state~$x(t)$ are known, %
then the auxiliary error terms~$e_i(t)$ may be expressed in terms of the state, the reference signal, the funnel functions and the derivatives of these. In the case of relative degree $r=2$ this means
\begin{equation}\label{eq:ei-Ei}
\begin{aligned}
	e_0(t) & = E_0(t,x(t)) = h(x(t)) - y_{\operatorname{ref}}(t), \\
	e_1(t) & = E_1(t,x(t)) = (L_f h)(x(t)) - \dot{y}_{\operatorname{ref}}(t)\\
&\qquad\qquad\qquad\quad + K_0(t,x(t)) E_0(t,x(t))
\end{aligned}
\end{equation}
with $K_0(t,x(t)) = 1/(1-\varphi_0^2(t) \|  E_0(t,x(t)) \|^2)$. Iteratively we may rewrite $e_i(t) =  E_i(t,x(t))$ for $i \in \{0,1,\ldots,r-1\}$.

\subsection{Funnel-MPC}

The funnel controller motivates to consider an associated Optimal Control Problem (OCP) with output constraints:
\begin{align}
	\underset{u \in \mathcal{L}^\infty([\hat{t},\hat{t}+T],\mathbb{R}^{p})}{\text{minimize}} \quad & \int_{\hat{t}}^{\hat{t}+T} \ell(t,x(t),u(t))\,\mathrm{d}t \nonumber \\
	\text{subject to } & \eqref{eq:Sys}, x(\hat{t}) = \hat{x} \text{ on $t \in [\hat{t},\hat{t}+T]$ and} \label{eq:OCP}\\
	\forall\,i \in \{0,\ldots,r-1\}\	 & \forall\, t \in [\hat{t},\hat{t}+T]:\ \| e_i(t) \| \le \varphi_i(t)^{-1}\nonumber
\end{align}
Next, we show existence of a feasible solution, for which the inequality constraint is inactive. To this end, we employ the funnel controller.

\begin{lem}\label{lem:funnel}
	Consider system~\eqref{eq:Sys} with relative degree~$r\in\N$, involutive distribution $\im g(x)$ and sign-definite high-gain matrix~$\Gamma(\cdot)$ as in~\eqref{eq:Gamma}.
	Moreover, let a reference trajectory $y_{\operatorname{ref}} \in \mathcal{C}^r
	([\hat{t},\hat{t}+T],\mathbb{R}^p)$ and positive funnel functions $\varphi_k \in \mathcal{C}^{r-k}
	([\hat{t},\hat{t}+T],\mathbb{R})$, $k = 0,\ldots, r-1$, be given. %
	In addition, choose the initial condition $x(\hat{t}) = \hat{x}$ such that $e_0,\ldots, e_{r-1}$ in~\eqref{eq:fun-con} satisfy
	\begin{equation}
		\forall\, i \in \{0,1,\ldots,r-1\}:\ \varphi_i(\hat t) \| e_i(\hat{t}) \| < 1. \label{eq:initialFeasibility}
	\end{equation}
	Then the funnel controller~\eqref{eq:fun-con} applied to~\eqref{eq:Sys} yields an initial-value problem %
	which has a unique solution~$x_{{\operatorname{FC}}}$ on~$[\hat{t},\hat{t}+T]$, %
	which is also an element of the feasible set of the OCP~\eqref{eq:OCP}. %
	Furthermore, there exists a unique~$\Theta >0$ such that
	\begin{equation}
		\forall\, i\!=\!0,\ldots,r\!-\!1:\ \|e_i(\hat{t}) \| + \Theta \leq \varphi_i(\hat{t})^{-1} \label{eq:errorboundTheta}
	\end{equation}
	and~\eqref{eq:errorboundTheta} is active, %
	i.e., there exist $j\in\{0,\ldots,r-1\}$ 
	such that~\eqref{eq:errorboundTheta} holds with equality for $i=j$. 
\end{lem}
\begin{pf}
	The assertion can be proved similar to~\cite[Thm.~3.1]{BergLe18}. We stress that compared to the latter work we do not require the BIBO property of the internal dynamics, because we only consider solutions on the compact time interval $[\hat{t},\hat{t}+T]$. More precisly, we may w.l.o.g.\ assume that the system is in Byrnes-Isidori form~\eqref{eq:BIF}, then it is clear that the closed-loop system~\eqref{eq:BIF},~\eqref{eq:fun-con} with initial value $x(\hat{t}) = \hat{x}$ has a solution on some interval~$[\hat t, \omega)$ such that $\|e_i(t)\|<\varphi_i(t)^{-1}$ for all $t\in[\hat t,\omega)$. As a consequence, $\xi$ is bounded on~$[\hat t, \omega)$ and integrating the second equation in~\eqref{eq:BIF} and using that~$q(\cdot,\cdot)$ is Lipschitz on~$[\hat t, \hat T]$ with Lipschitz bound $L>0$ (where $\hat T>\omega$ is arbitrary) we obtain
\begin{align*}
    \|\eta(t)\| &\le \|\eta(\hat t)\| + \int_{\hat t}^t \|q(0,0)\| + L \|(\xi(s),\eta(s))\| {\rm d}s \\
    &\le \|\eta(\hat t)\| + M_1 (\omega - \hat t) + M_2 \int_{\hat t}^t \|\eta(s)\| {\rm d}s
\end{align*}
for all $t\in[\hat t,\omega)$ and some $M_1, M_2>0$. Then Gronwall's lemma implies that~$\eta$ is bounded on~$[\hat t,\omega)$ and hence, analogous to~\cite[Thm.~3.1]{BergLe18}, we may show that the solution can be extended to~$[\hat t, \hat T]$. Since $\hat T$ was arbitrary we may choose $\hat T = \hat t + T$. Furthermore, because all involved functions in~\eqref{eq:Sys} are sufficiently smooth, they are locally Lipschitz and hence the solution~$(\xi,\eta)$ is unique. This guarantees uniqueness of~$\Theta$.\hspace*{31mm} \qed
\end{pf}

Based on Lemma~\ref{lem:funnel}, we define the map $\Psi: \cD\times \R_{\ge 0} \rightarrow \mathbb{R}_{>0}$ with $\Psi(\hat t, \hat x, T) = \Theta$, where $\cD\subseteq \R_{\ge 0}\times\R^n$ is the set of points $(\hat t,\hat x)$ which satisfy~\eqref{eq:initialFeasibility}. Hence, $\Psi$ is well-defined, provided the assumptions of Lemma~\ref{lem:funnel} hold. %
Note that $\Theta$ may have been chosen to be an element of $\mathbb{R}^r$ using the same line of reasoning applied componentwise.\\
We exploit the map~$\Psi$ and, thus, indirectly the funnel controller~\eqref{eq:fun-con}, in order to ensure initial and recursive feasibility of the following MPC scheme.

\noindent\textbf{Algorithm: Funnel-MPC}\\
\noindent\textbf{Given}: system \eqref{eq:Sys} with relative degree~$r$, reference signal $y_{\operatorname{ref}} \in \mathcal{C}^r(\R_{\geq 0},\mathbb{R}^p)$ and funnel functions $\varphi_k \in \mathcal{C}^{r-k}(\mathbb{R}_{\geq 0},\mathbb{R})$, $k \in \{0,1,\ldots,r-1\}$.\\
\noindent\textbf{Set} the time shift $\delta >0$, the prediction horizon $T \ge \delta$ and the current time $\hat{t} := 0$.
\begin{enumerate}
	\item[(a)] Set $\hat{x} := x(\hat{t})$
	\item[(b)] Solve the OCP~\eqref{eq:OCP} subject to the additional constraint
		\begin{equation}
			\forall\, i=0,\ldots,r-1:\ \| e_i(\hat{t}+\delta) \| \leq \varphi_i(\hat t + \delta)^{-1} - \Psi(\hat{t},\hat{x},T) \label{eq:feasibilityConstraint}
		\end{equation}
	 in order to compute an (approximately) optimal control function $u^\star \in \mathcal{L}^\infty([\hat{t},\hat{t}+T],\mathbb{R}^p)$.
	 \item[(c)] Implement the feedback law $\mu: [\hat{t},\hat{t}+\delta) \times \mathbb{R}^n \rightarrow \mathbb{R}^p$, $\mu(t,\hat{x}) := u^\star(t)$ at system~\eqref{eq:Sys}, increase $\hat{t}$ by $\delta$ and go to step~(a)
\end{enumerate}

Next, we show initial and recursive feasibility of the proposed MPC scheme.
\begin{thm}\label{thm:feasibility}
	Consider system~\eqref{eq:Sys} with relative degree~$r\in\N$, involutive distribution $\im g(x)$ and sign-definite high-gain matrix~$\Gamma(\cdot)$ as in~\eqref{eq:Gamma}. %
	Further, let a reference trajectory $y_{\operatorname{ref}} \in \mathcal{C}^r
	(\R_{\ge 0},\mathbb{R}^p)$ and positive funnel functions $\varphi_k \in \mathcal{C}^{r-k}
	(\R_{\ge 0},\mathbb{R})$, $k \in \{0,1,\ldots,r-1\}$, be given. %
	If the initial data $(\hat t,\hat{x}) = (0,x(0))$ satisfies condition~\eqref{eq:initialFeasibility}, then Funnel-MPC is initially and recursively feasible, i.e., %
	the feasible set of the OCP~\eqref{eq:OCP} augmented by the feasibility constraint~\eqref{eq:feasibilityConstraint} is non-empty at time $\hat{t} = 0$ %
	and at each successor time $\hat{t} \in \delta \mathbb{N}$.
\end{thm}
\begin{pf}
	The proof is a straightforward consequence of the proposed construction. %
	Condition~\eqref{eq:feasibilityConstraint} ensures that the funnel controller yields %
	a feasible solution of the OCP~\eqref{eq:OCP} augmented by the feasibility constraint~\eqref{eq:feasibilityConstraint} in view of Lemma~\ref{lem:funnel}. %
	Since this constraint is incorporated in the OCP~\eqref{eq:OCP} and the feasible set at time $t = \hat{t} + \delta$ is closed, %
	also the optimal control~$u^\star$ satisfies~\eqref{eq:feasibilityConstraint}. %
	This, however, implies that condition~\eqref{eq:initialFeasibility} holds at the successor time. %
	Since this line of reasoning can also be applied at each subsequent time $\hat{t} \in \delta \mathbb{N}$, recursive feasibility follows. \hspace*{19mm}\qed
\end{pf}

We stress that in Theorem~\ref{thm:feasibility}, in order to apply the Funnel-MPC scheme, we do not require that $y_{\rf}$, $\varphi_k$ and its derivatives are bounded, which is needed for the application of the funnel controller~\eqref{eq:fun-con}. Theorem~\ref{thm:feasibility} ensures well-posedness of the MPC closed loop system, provided that~\eqref{eq:initialFeasibility} holds, i.e., %
the initial state is contained in the interior of the funnel. %
The funnel functions~$\varphi_k$ can always be chosen such that this holds. Although not explicitly allowed in Theorem~\ref{thm:feasibility}, it is also possible to use ``infinite funnels'', i.e., $\varphi_k(0)=0$, so that any initial value is feasible.

\subsection{Stage Cost}

The stage costs~$\ell(\cdot,\cdot,\cdot)$ in~\eqref{eq:OCP} have not been specified so far and we like to use costs, which only depend on the auxiliary errors and the control effort, e.g., using the rewritten auxiliary errors in~\eqref{eq:ei-Ei},
\begin{align}
	\ell(t,x,u) := \sum_{i=0}^{r-1} \| E_i(t,x) \|^2 + \lambda \| u \|^2 \label{eq:stageCostClassical}
\end{align}
with regularization parameter $\lambda >0$. %
However, the stage costs~\eqref{eq:stageCostClassical} penalize the errors independently of the specified funnel. %
If e.g.\ a batch process is considered, then a wider funnel corresponds to less emphasis on the error. %
Hence, penalizing the distance of the error from the funnel boundary seems a more reasonable approach. %
A straightforward way to achieve this is to use the funnel gains~$k_i$ as in~\eqref{eq:fun-con} (which have to be computed for evaluating the error anyway), i.e.,
\begin{align}
	\ell(t,x,u) := \sum_{i=0}^{r-1} \frac {1}{1-\varphi_i(t)^2\| E_i(t,x) \|^2} + \lambda \| u \|^2. \label{eq:stageCostFunnel}
\end{align}
One alternative would be to focus on $e_0(\cdot,\cdot)$ and neglecting \textit{higher order} terms. Moreover, we emphasize that both performance measures can be evaluated without knowing the state --- only the output, the control, and the current time are needed.

\section{Funnel-MPC: A Numerical Case Study}

In this section, we illustrate the Funnel-MPC scheme by a numerical case study consisting of three parts: First, we investigate funnel control and its implementation with zero-order hold (ZOH), i.e.,
\begin{equation*}
	u_{\operatorname{ZOH}}(t) = u_{\operatorname{FC}}( \lfloor t /\tau \rfloor \tau ) 
\end{equation*}
meaning that the control is only updated every $\tau$ time units, which leads to a sampled-data system with ZOH. %
Then we investigate the Funnel-MPC scheme with time shift~$\delta = \tau$ and show that it outperforms the controller~$u_{\operatorname{ZOH}}$ w.r.t.\ required sampling rate and performance. %
Finally, we present a combination of the funnel controller and Funnel-MPC, which is applicable if the model is not available, but its parameters must first be identified.\\
Throughout this section, we consider an example of a mass-spring system mounted on a car from~\cite{SeifBlaj13}, see Fig.~\ref{fig:FunnelZOH}. The equations of motion (where the control input is the force acting on it) are given by
		\begin{align*}
			\begin{bmatrix}
				m_1+m_2 & m_2 \cos(\alpha) \\
				m_2 \cos(\alpha) & m_2
			\end{bmatrix} \hspace*{-0.1cm} \begin{pmatrix}
				\ddot{x}(t) \\
				\ddot{s}(t)
			\end{pmatrix} \!+\! \begin{pmatrix}
				0 \\
				ks(t) + d\dot{s}(t)
			\end{pmatrix} \!=\! \begin{pmatrix}
				u(t) \\
				0
			\end{pmatrix}
		\end{align*}
and the output is chosen as the horizontal position of the mass on the ramp,
\[
  y(t) = x(t) + s(t) \cos(\alpha).
\]
Clearly, the mass on car system can be rewritten in the form~\eqref{eq:Sys} and, as shown in~\cite{BergLe18}, the system has relative degree $r=2$ for $0 < \alpha < \pi/2$ and $r=3$ for $\alpha = 0$; with positive high-gain matrix~$\Gamma$ in both cases. For the simulations we use the parameters $m_1 = 4$ and $m_2 = 1$ for the mass of the car and the mass moving on the ramp, resp., $k = 2$ and $d = 1$ for the coefficients of the spring and damper, resp., and the initial values $x(0) = \dot x(0) = s(0) = \dot s(0) = 0$.

\subsection{Funnel Control with ZOH Sampling}\label{Subsection:FunnelZOH}

\begin{figure}[!htb]
	\begin{center}
		\includegraphics[trim=2cm 4cm 5cm 15cm,clip=true,width=0.24\textwidth]{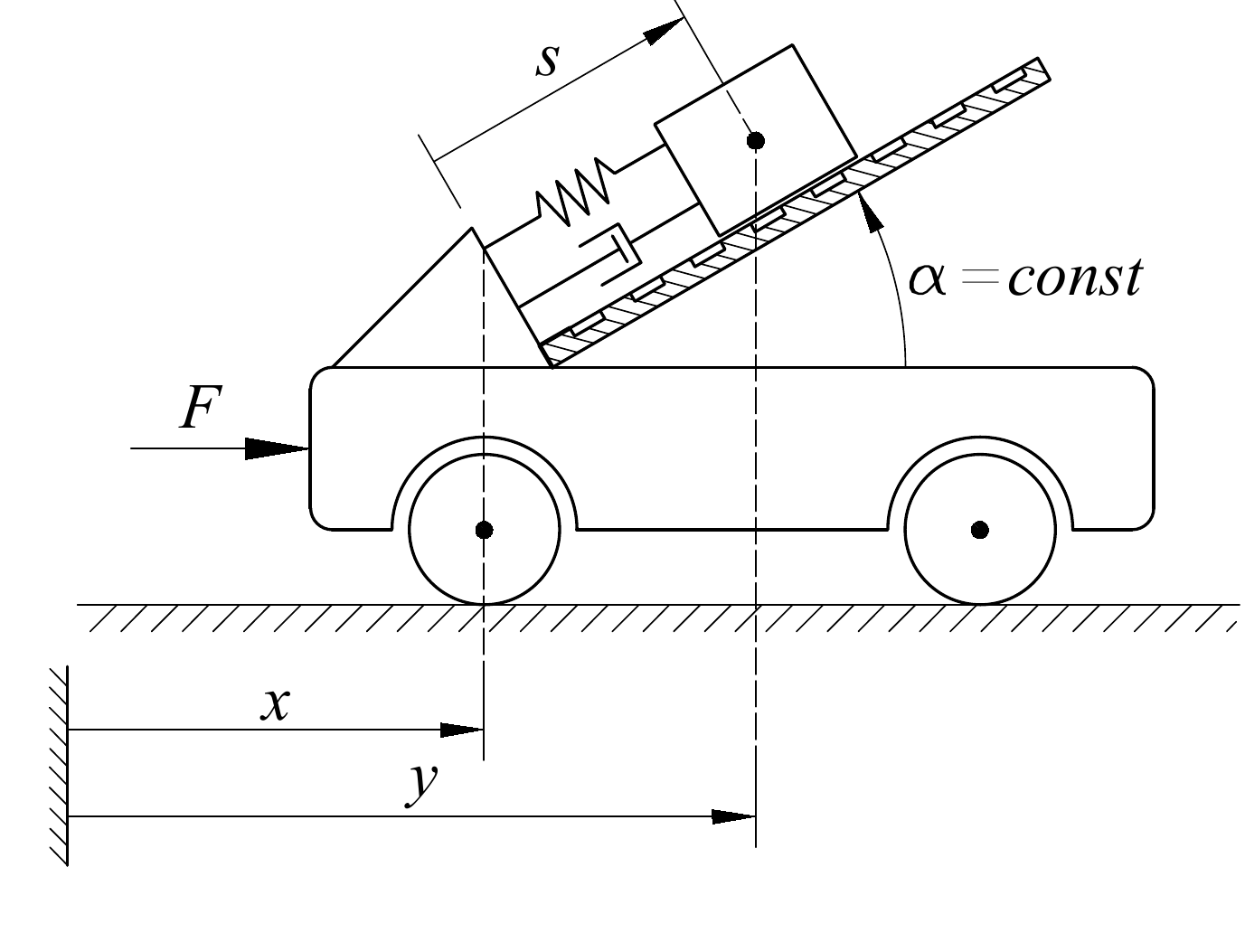}
		\includegraphics[width=0.24\textwidth]{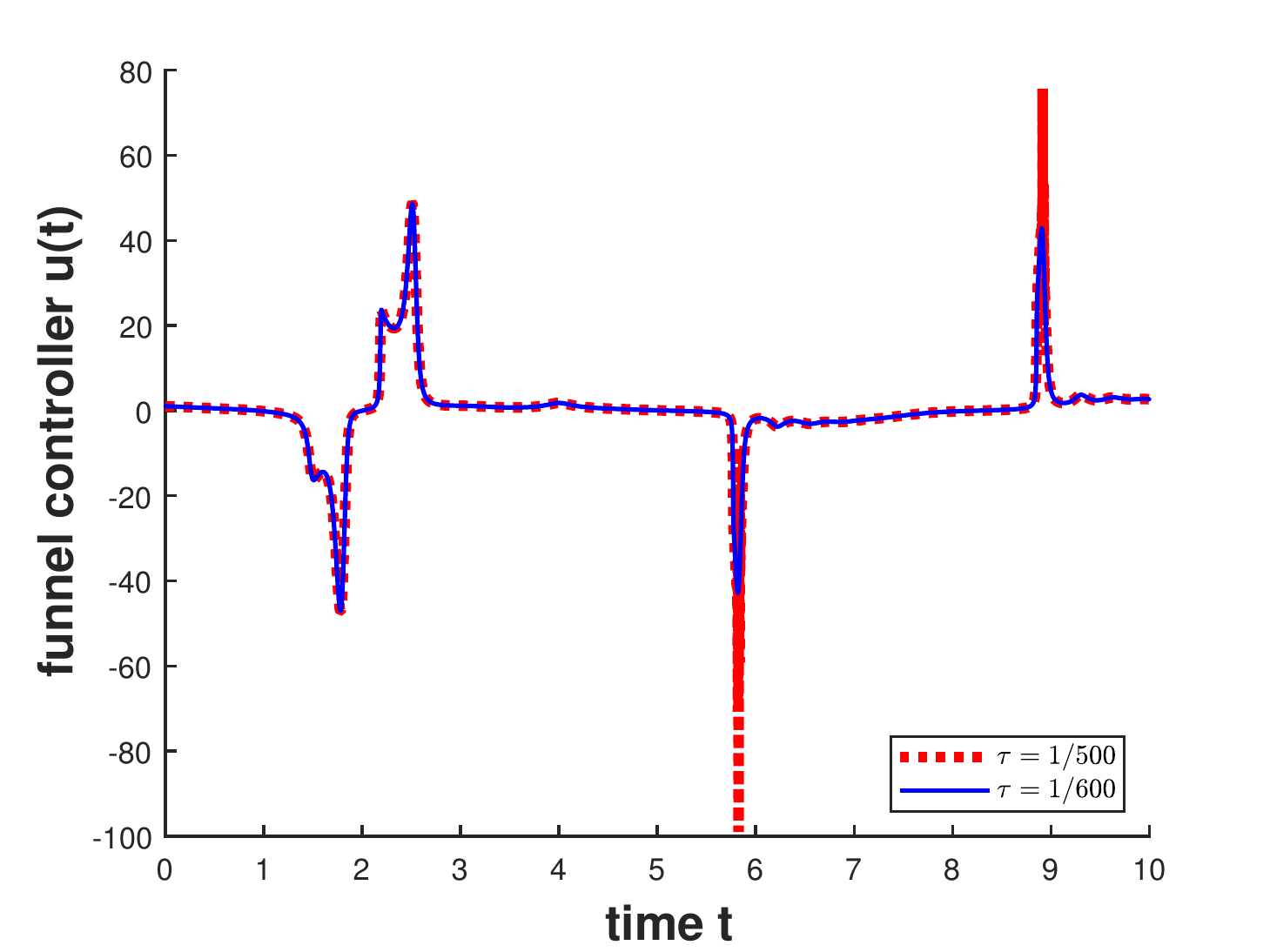}
		\includegraphics[width=0.24\textwidth]{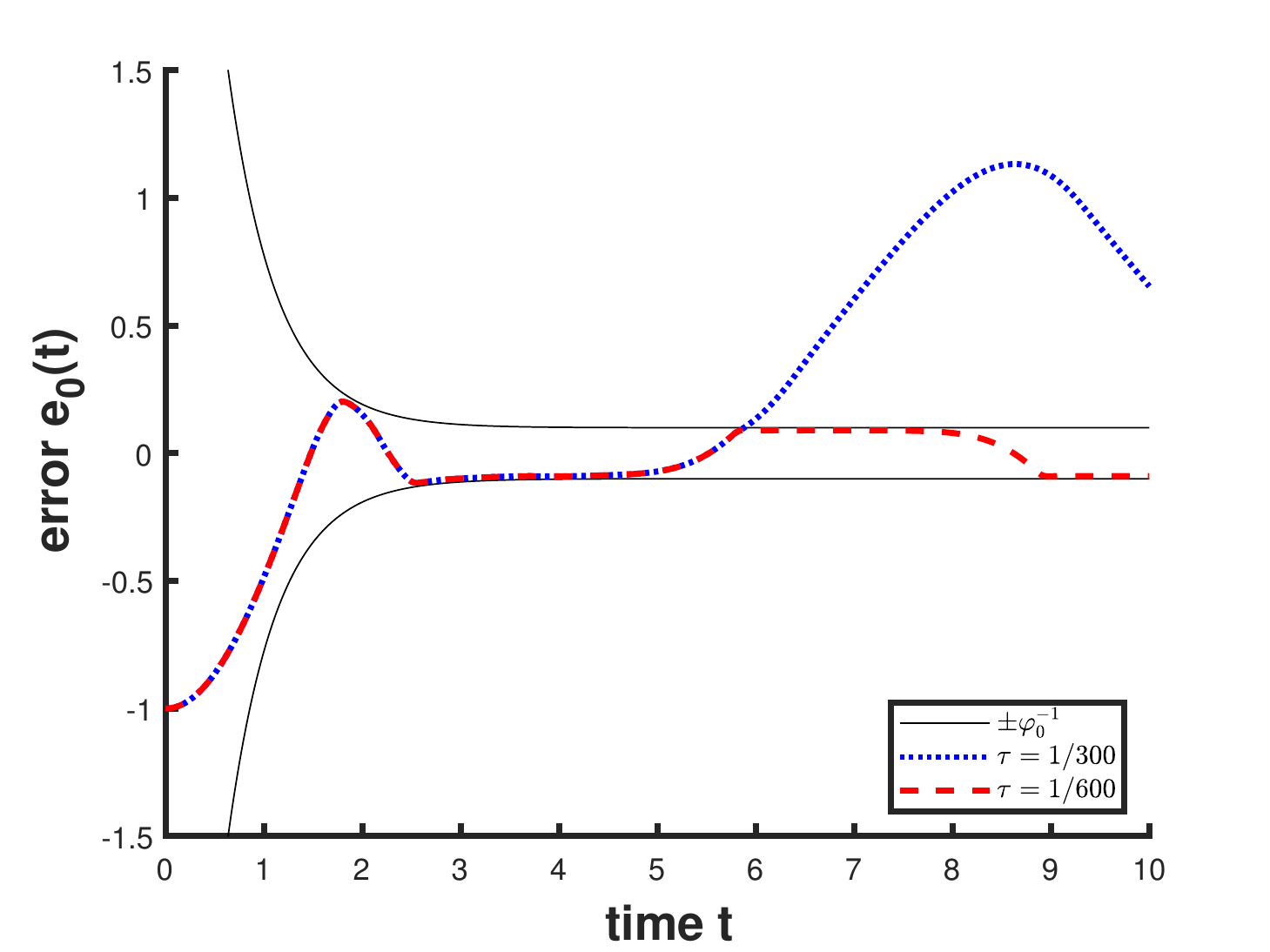} 
		\includegraphics[width=0.24\textwidth]{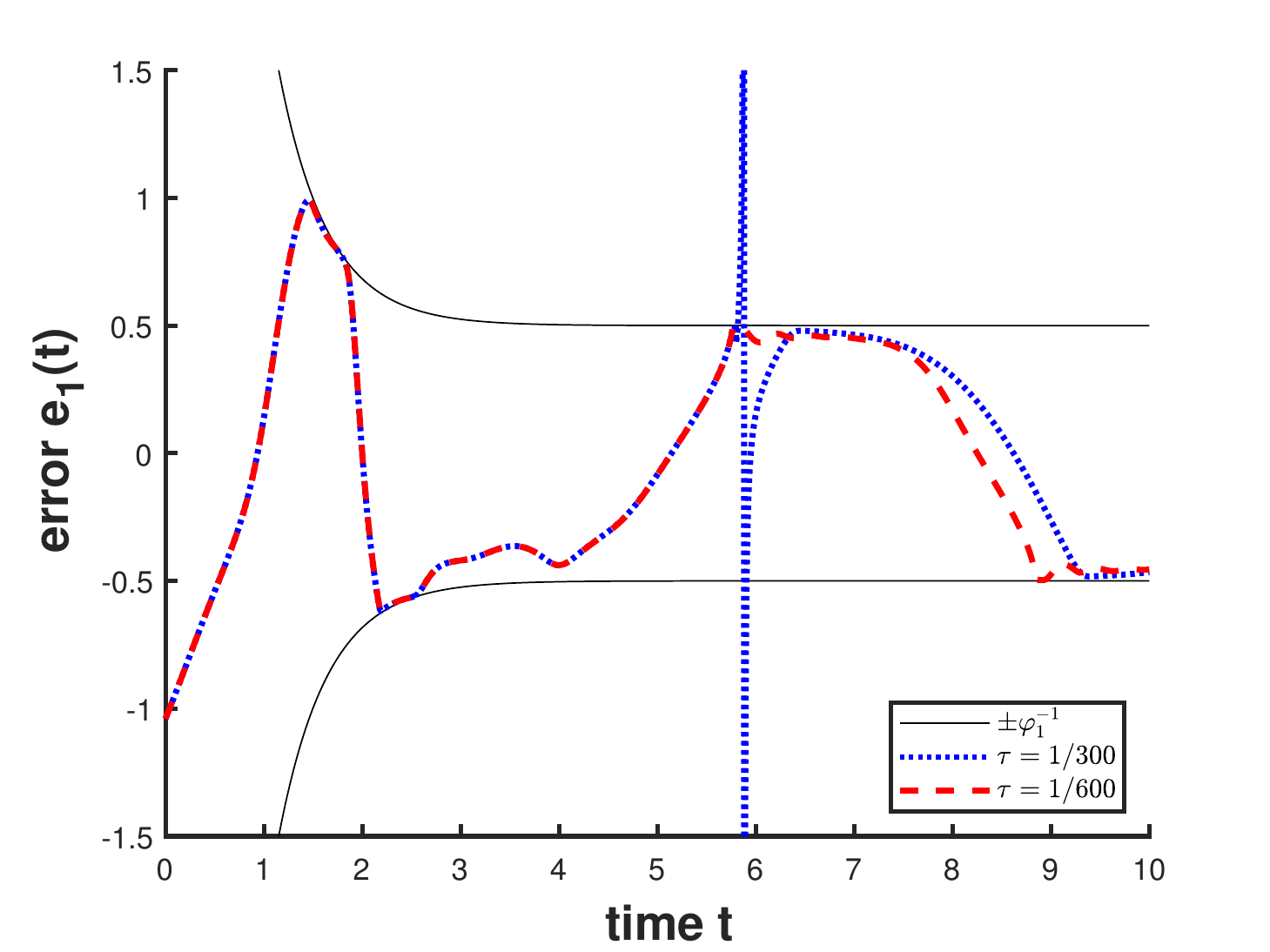}
		\caption{Schematic sketch of the mass on car system and ZOH implementation of the funnel controller: control effort (above) and output errors (below).}
		\label{fig:FunnelZOH}
	\end{center}
\end{figure}

In this subsection, we consider the ZOH implementation of the funnel controller. We neglect potential additional difficulties inferred from the requirement of computing derivatives of the output by providing the exact values at each sampling instant using the representation of the auxiliary errors in terms of the state as in~\eqref{eq:ei-Ei}.\\
The simulation on the time interval $[0,10]$ is performed using the \textsc{Matlab}-routine \texttt{ode45}. We first choose $\alpha = \pi/4 \in (0,\pi/2)$ (relative degree two) and the funnels
\begin{equation}\label{eq:funnelLinear}
\begin{aligned}
	\varphi_0(t) := (0.1 + 5 \exp(-2t))^{-1},  \\
	\varphi_1(t) := (0.5 + 10 \exp(-2t))^{-1}.
\end{aligned}
\end{equation}
As shown in Fig.~\ref{fig:FunnelZOH}, feasibility (i.e., error evolution within the funnel boundaries) is not maintained for a sampling rate $\tau = 1/300$. While feasibility is achieved for $\tau = 1/500$, the control signal is deteriorated as the range of the control values is significantly larger. For the sampling rate $\delta = 1/600$ the continuous-time performance (cf.~\cite{BergLe18}) is essentially recovered.
For relative degree $r = 3$ ($\alpha = 0$), $\varphi_0$ as in~\eqref{eq:funnelLinear} and $\varphi_1(t) = \varphi_2(t) = (0.05 + 1.4 \exp(-t))^{-1}$, 
feasibility is ensured for sampling rate $\delta = 1/700$. However, a shattering of the control can be observed, which results in oscillations within the time intervals $[5.6,6.7]$ and $[8.6,9.6]$. Here, a sampling rate of $\tau = 1/1200$ is needed in order to recover the performance of the funnel controller with ZOH. We remark that a redesign of the computed control signal for the digital implementation may be beneficial, see e.g.~\cite{GrunWort08,GrunWort08b} and~\cite{MonaNorm10}.


\subsection{Funnel-MPC Scheme}

\begin{figure}[!htb]
	\begin{center}
		\includegraphics[width=0.24\textwidth]{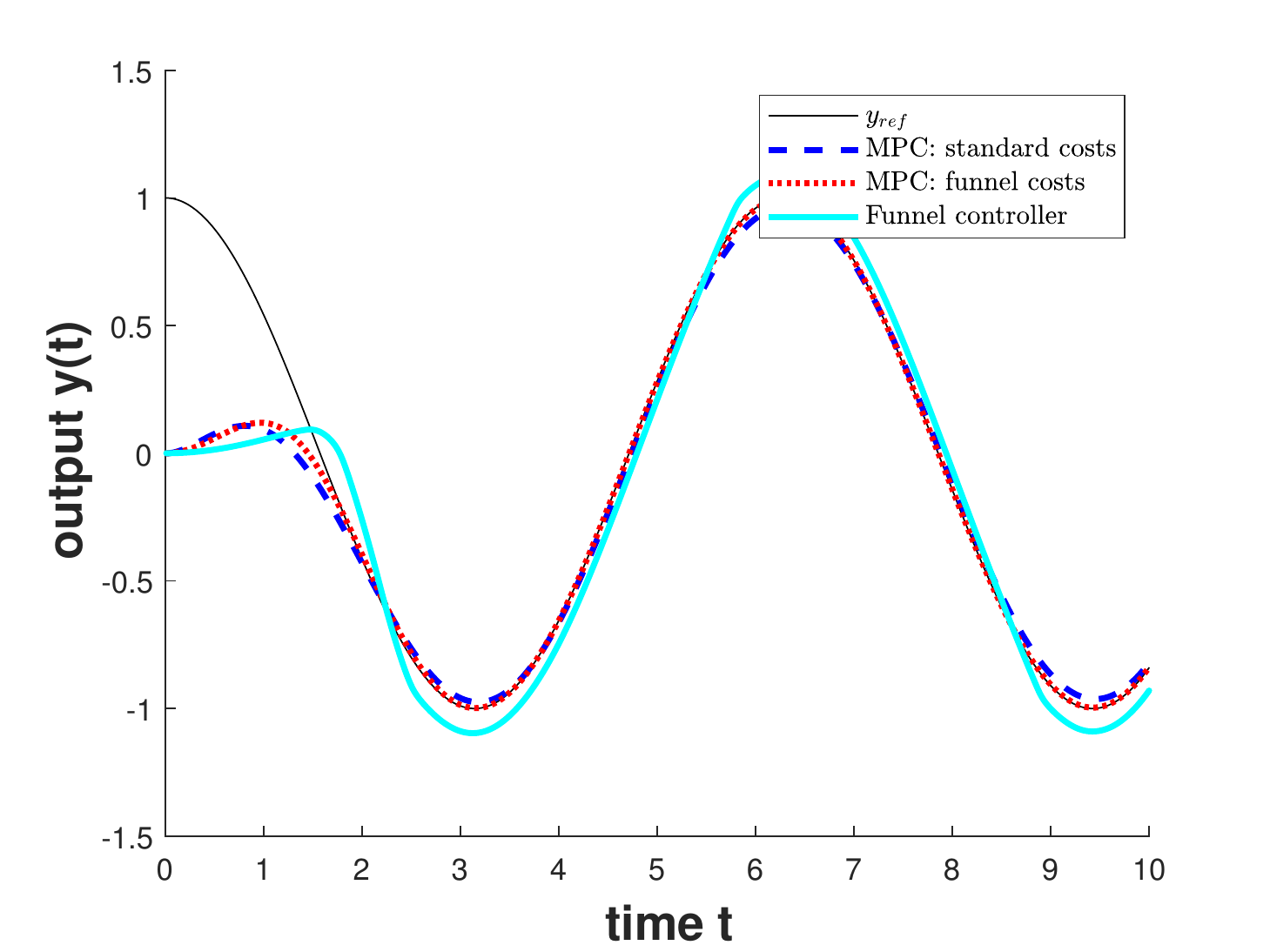}
		\includegraphics[width=0.24\textwidth]{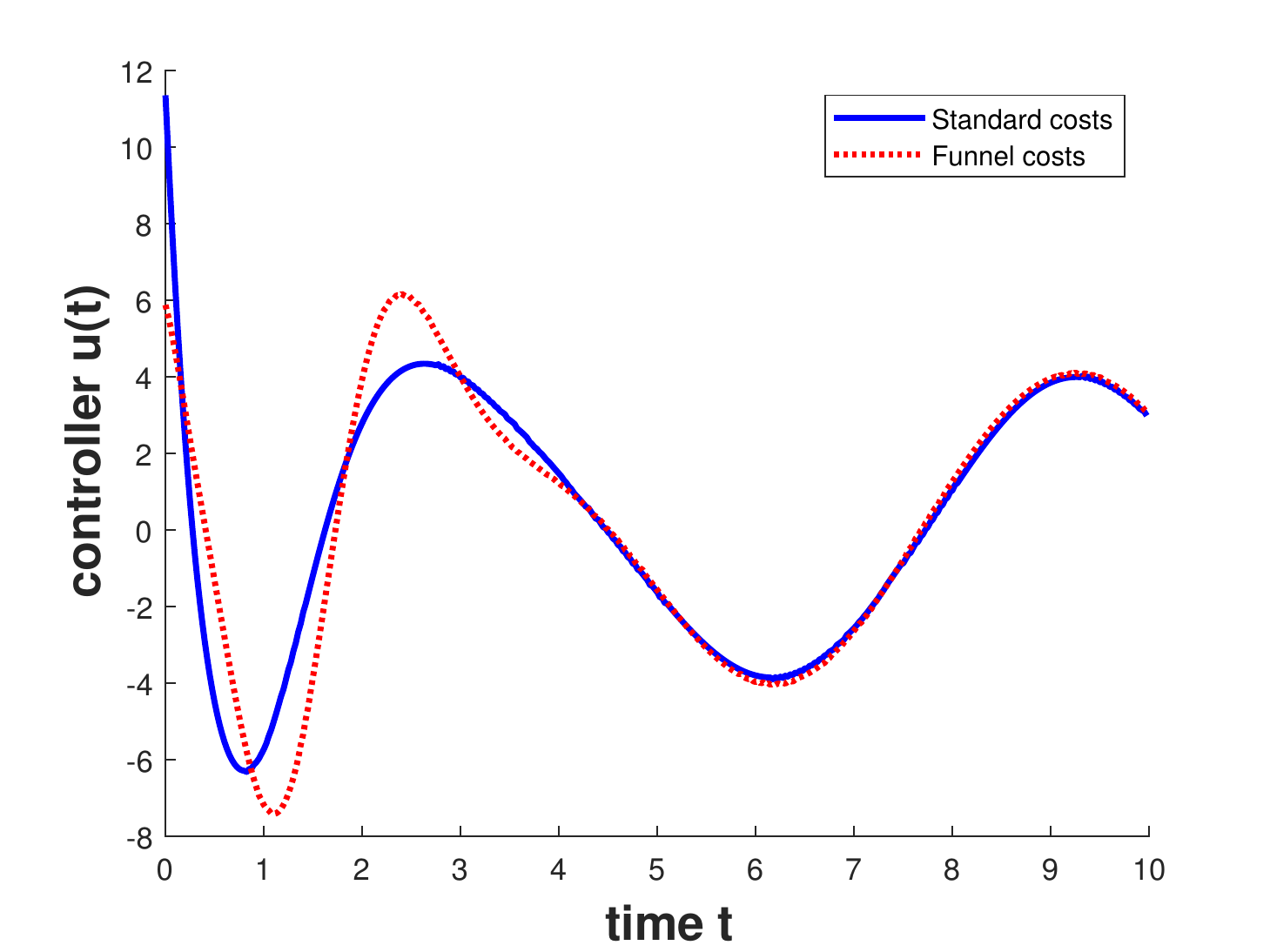}
		\includegraphics[width=0.24\textwidth]{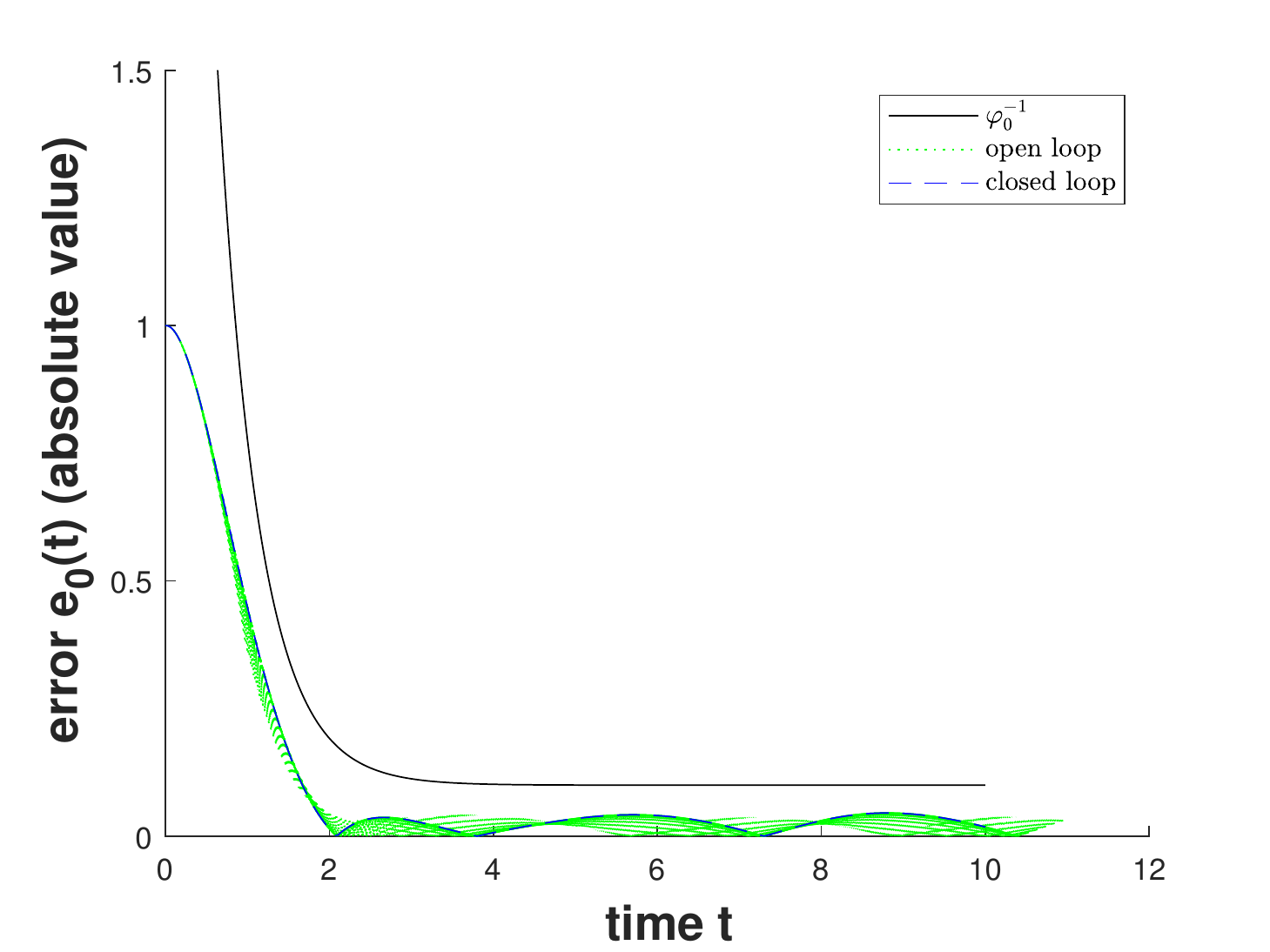}
		\includegraphics[width=0.24\textwidth]{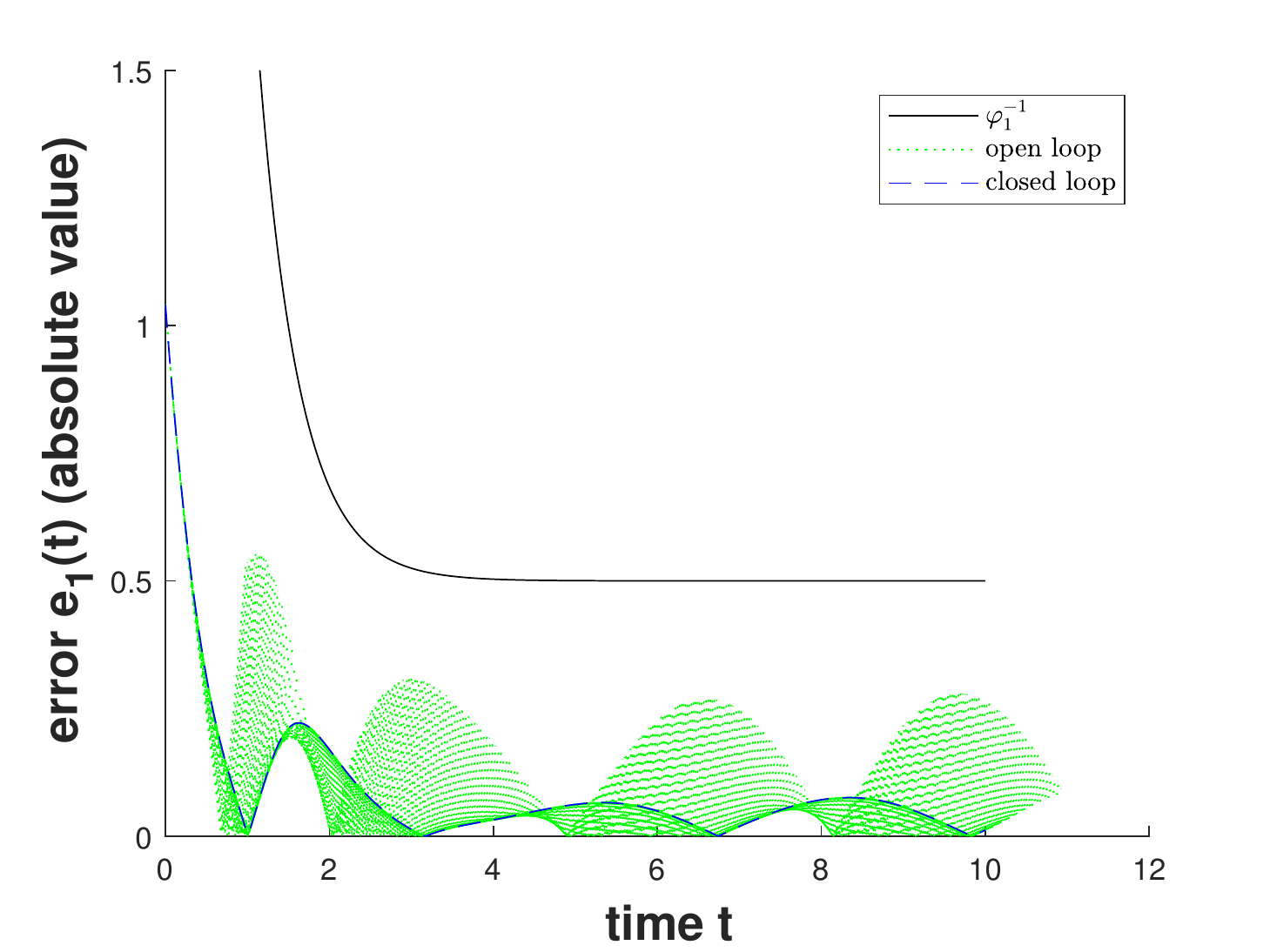}
		\includegraphics[width=0.24\textwidth]{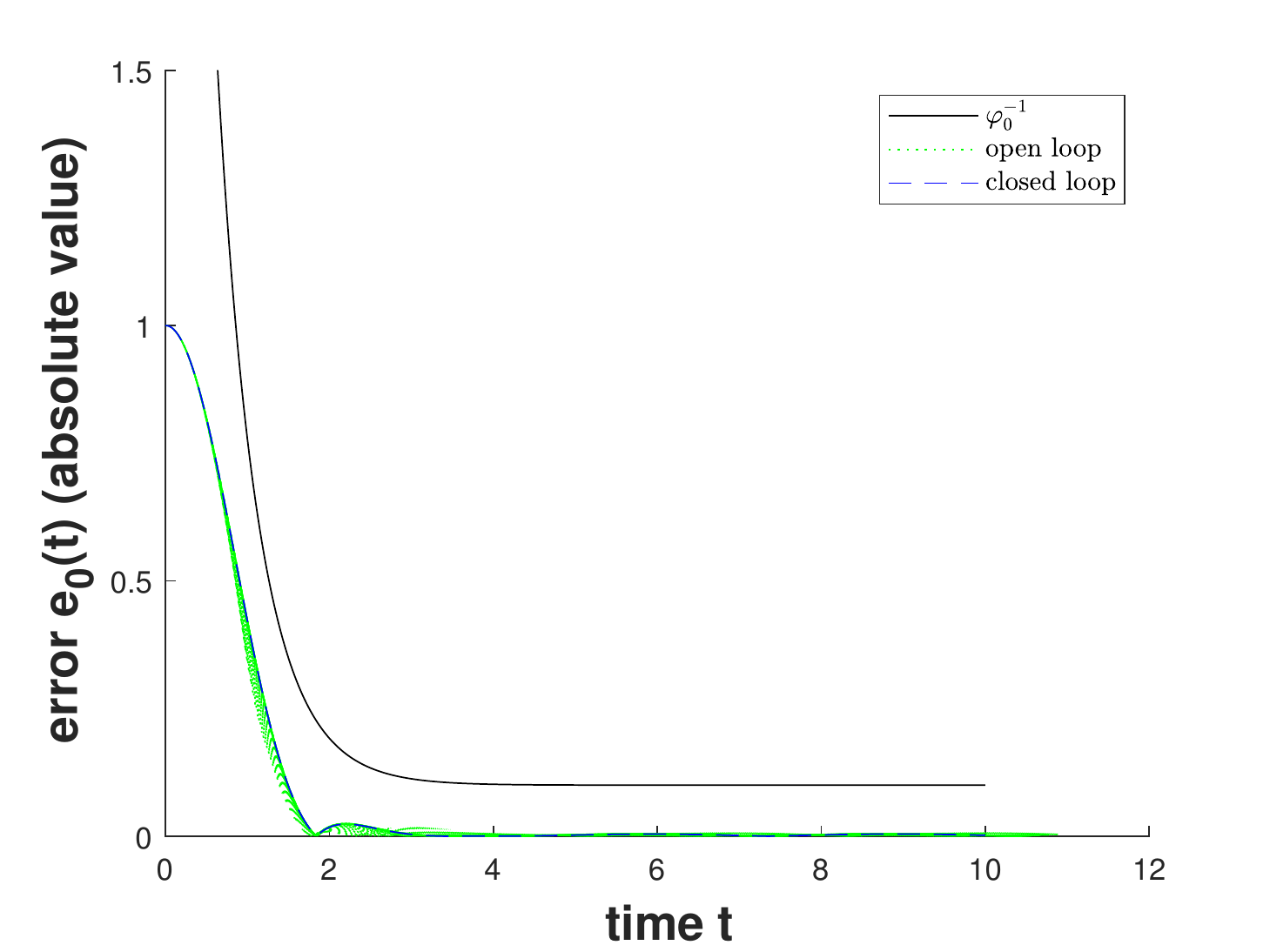}
		\includegraphics[width=0.24\textwidth]{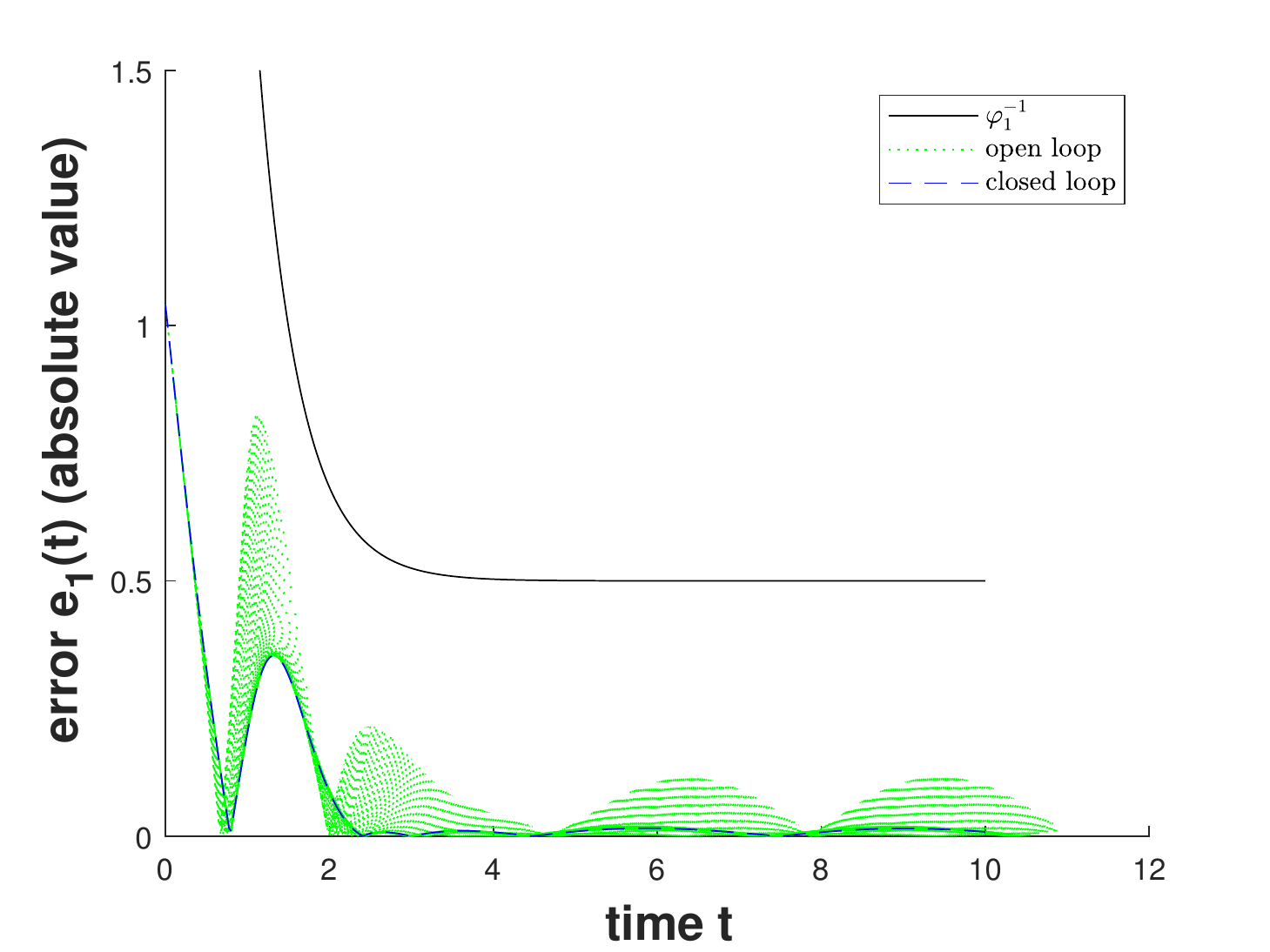}
		\caption{Output and control effort for the Funnel-MPC scheme applied to the mass on car system (top), %
		trajectories of the absolute value of the errors $e_0$, $e_1$ with stage costs~\eqref{eq:stageCostClassical} (middle) and~\eqref{eq:stageCostFunnel} (bottom).}
		\label{fig:MPC}
	\end{center}
\end{figure}

We use MPC without stabilizing terminal constraints or costs, see e.g.\ the textbook~\cite{GruePann17} for details. Moreover, we refer to~\cite{WortRebl14} and the references therein for a detailed analysis and discussion of the connection between continuous-time and discrete-time MPC schemes without stabilizing terminal constraints and costs. Sampled-data systems with ZOH are in particular treated in~\cite{WortRebl15}.\\
Again, we consider the mass on car system with $\alpha = \pi/4$ and, thus, relative degree $r=2$. For the Funnel-MPC scheme we use the prediction horizon $T=N\delta$ with $N=41$ and time shift $\delta = \tau = 1/40$. First, we employ the ``standard'' stage costs~\eqref{eq:stageCostClassical} with weighting factor $\lambda = 0.005$. Note that the sampling rate is much lower (factor~$15$) compared to the ZOH implementation of the funnel controller. Hence, the funnel controller is able to adjust its control signal significantly more often. Nevertheless, MPC yields a feasible control input and a drastically increased performance measure
\begin{align*}
	\sum\nolimits_{i=0}^{10/\delta} \ell(i\tau,x_{\operatorname{MPC}}(i\tau),\mu_{MPC}(i\tau,x_{\operatorname{MPC}}(i\tau)))
\end{align*}
of~$1.5435$ 
on the simulated time interval, i.e., less than $22\%$ of the costs evaluated along the funnel control and output error trajectories, which yield aggregated costs of $7.0900$. This can also be observed from Fig.~\ref{fig:MPC}: the output is tracked more accurately and the range of employed control values is much smaller, i.e., approximately $[-6,12]$ compared to $[-50,50]$ for the funnel controller. Here, both the Euler and the \textsc{Matlab}-routine \texttt{ode45} yield essentially the same results for the MPC closed loop in our numerical investigations (Euler yields aggregated costs of~$1.5511$ instead of~$1.5435$). Hence, we present the results computed with Euler in the following.
Next, we consider the Funnel-MPC scheme with the 
stage costs~\eqref{eq:stageCostFunnel}, which penalize the distance of the auxiliary errors to the funnel boundaries instead of their absolute values (we use the term \textit{funnel costs} in Fig.~\ref{fig:MPC}). It can be observed that already the open-loop predictions provide a larger distance to the funnel boundary (the performance measure decreases from $107.8013$ to $21.0963$). %
This controller design is particularly interesting since it may allow to eschew the funnel constraints and the additional feasibility constraint~\eqref{eq:feasibilityConstraint} in the OCP~\eqref{eq:OCP}. We conjecture that feasibility of the Funnel-MPC scheme without incorporating the funnel constraints and the feasibility constraint~\eqref{eq:feasibilityConstraint} can be proved by using optimality of the predicted open-loop trajectories.\\
For relative degree $r=3$ ($\alpha = 0$), %
feasibility along the MPC closed-loop trajectory is preserved for the sampling rate $\delta = \tau = 1/40$ %
(factor~$30$ compared to the funnel controller). %
However, the improvement of the performance measure is smaller, i.e.,~$1.8171$ instead of~$2.4389$. %
\begin{rem}
If a model of the system is available, then a reference control signal~$u_{\operatorname{ref}}$ which exactly reproduces the reference~$y_{\rf}$ can be obtained using system inversion techniques and a corresponding feedforward-control input can be computed. However, if the initial error $e(0)$ is nonzero, then this approach may lead to large errors and must be combined with feedback control techniques. Such a combination with the funnel controller has been proposed in~\cite{BergOtto19} and was simulated for the mass on car system. In a corresponding combination with a MPC scheme the difference $u(t) - u_{\rf}(t)$ could be penalized instead of~$u(t)$. However, even~$u_{\operatorname{ref}}(t)$ takes values in the interval $[-8.5,4.5]$ with a peak at $t = 0$. %
	This explains the high value at the beginning caused by the large initial error due to the choice of the initial value. %
	Hence, the proposed Funnel-MPC controller yields a satisfactory range of control values.
\end{rem}
%
While the above numerical experiments look promising, one should keep in mind that the computational effort of MPC is much higher compared to funnel control 
and the controller is non-causal, i.e., predictions depending on future inputs have to be computed (which may be done using a system model). The latter drawback motivates the investigations in the following subsection.

\subsection{Model identification during runtime}

\begin{figure}[!htb]
	\begin{center}
		\includegraphics[width=0.24\textwidth]{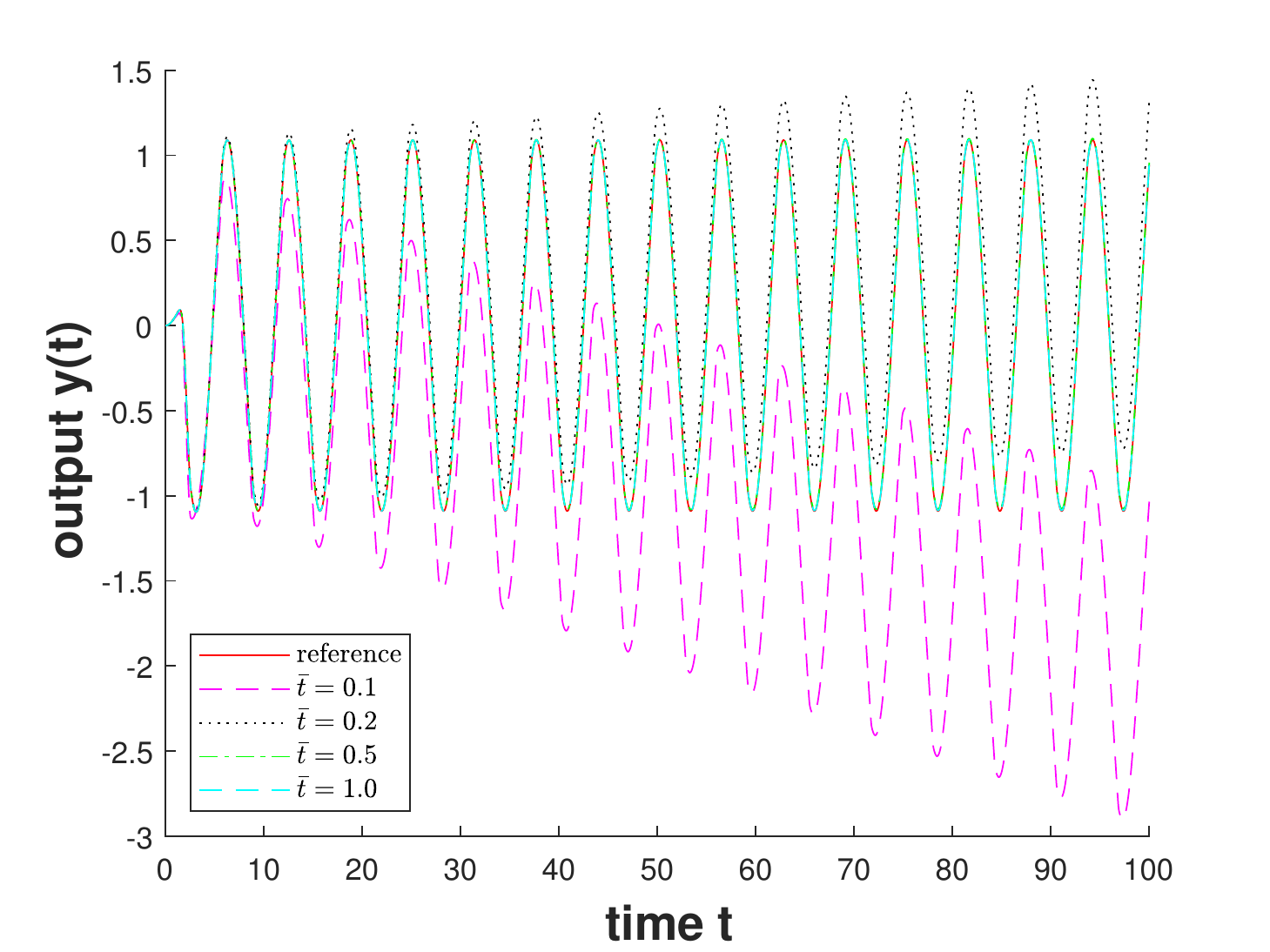}
		\includegraphics[width=0.24\textwidth]{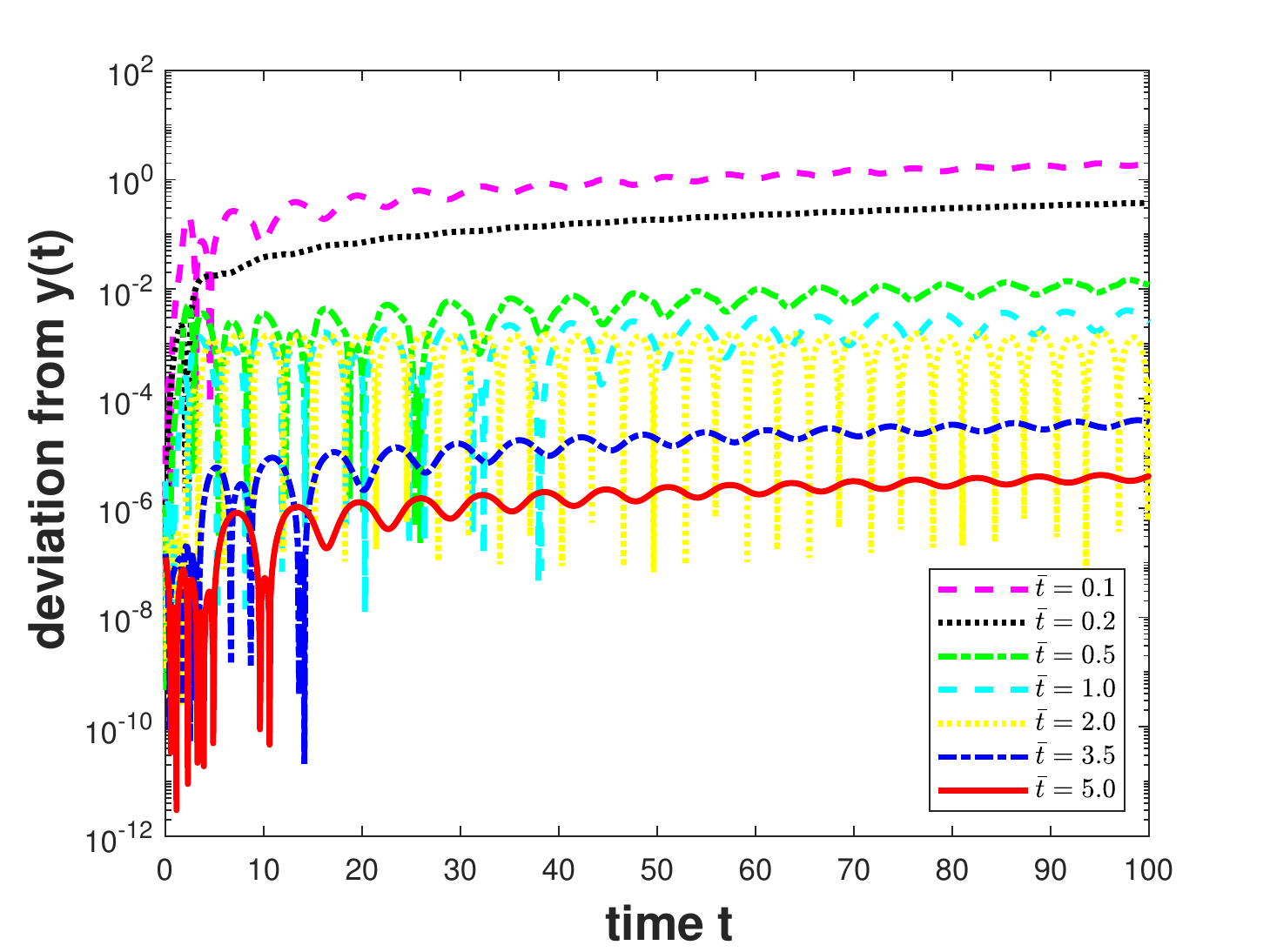}
		\caption{Open-loop trajectories based on the identified system parameters $\alpha$, $m_1$, $m_2$, $k$, $d$ and $z^0$ after $\bar{t}$ time units.}
		\label{fig:predictions}
	\end{center}
\end{figure}

The assumption that the system dynamics and the initial value are known restrict the applicability of Funnel-MPC, especially when compared to the model-free funnel controller. In this subsection, we present a methodology to resolve this drawback by using funnel control during a learning phase, where the system model and the initial/current state are identified, and Funnel-MPC as soon as the model is sufficiently well known.\\
We assume knowledge of the structure of the mass on car system, but only limited information on the parameters, i.e., $\alpha \in [0,\pi/2]$, $m_1 \in [2,6]$, $m_2 \in [0.5,1.5]$, $k \in [1,3]$, $d \in [0.5,1.5]$, and the initial value $z^0 = (x(0),\dot x(0), s(0), \dot s(0))^\top$, i.e.,
\begin{align*}
	z^0 \in [-2.5,3.5] \times [-1,1] \times [-2.75,3.25] \times [-1,1].
\end{align*}
As it can be inferred from the simulations in Subsection~\ref{Subsection:FunnelZOH}, the funnel controller yields a satisfactory behaviour on the time interval~$[0,1]$. Hence, we apply the funnel controller with step size $\tau = 10^{-3}$ in order to determine an input-output vector $(u_{\operatorname{FC}},y_{\operatorname{FC}}) = (u(i\tau),y(i\tau))_{i=0}^{100q}$, $q \in \{1,2,5,10\}$, and solve the minimization problem
\begin{align*}
	\underset{\alpha,m_1,m_2,k,d,z^0}{\text{minimize}} \quad & \| \hat{y} - y_{\operatorname{FC}} \|^2 \nonumber \\
	\text{subject to } & z(0) = z^0 \text{ and for all $i \in \{1,2,\ldots,100q\}$} \\
		 & z(i) = \varphi(\tau;z(i-1),u_{\operatorname{FC}}(i-1))
		 \nonumber\\
		 & y(i) = [1, \cos(\alpha), 0, 0] z(i) \nonumber
\end{align*}
where $z=(x,\dot x, s,\dot s)^\top$ denotes the state of the mass on car system and $\varphi(\tau;z(i-1),u_{\operatorname{FC}}(i-1))$ denotes the solution of it with initial value~$z(i-1)$ and constant control $u(t) \equiv u_{\operatorname{FC}}(i-1)$ after $\tau$\,time units. We obtain the quadratic and maximal error of the predicted trajectories on the time interval~$[0,100]$ shown in Table~\ref{tb:predictions}.\\
\begin{table}[!htb]
	\begin{center}
		\begin{tabular}{|c|c|c|c|}\hline
			& $\| y - \hat{y} \|_2$ & $\| y - \hat{y} \|_\infty$ \\\hline
			$\bar{t} = 0.1$ & 358.5043074855 &	 1.9866952202 \\
			$\bar{t} = 0.2$ &  68.5598064742 &	 0.3769510455 \\
			$\bar{t} = 0.5$ &   2.2736367216 &	 0.0147993864 \\
			$\bar{t} = 1.0$ &   0.6038075070 &	 0.0040749474 \\
			$\bar{t} = 2.0$ &   0.3149716888 &	 0.0015685639 \\
			$\bar{t} = 3.5$ &   0.0066075315 &	 0.0000401394 \\
			$\bar{t} = 5.0$ &   0.0006860503 &	 0.0000039939 \\ \hline
		\end{tabular}
	\end{center}
    \ \\[-2mm]
	\caption{Prediction errors on the time interval $[0,100]$ after \textit{learning} on the interval $[0,\bar{t}]$.}
	\label{tb:predictions}
\end{table}
Clearly, on the prediction interval~$[\hat{t},\hat{t}+T]$ of length one, the error is significantly smaller, as it can be inferred from the open-loop error trajectories, cf.\ Fig.~\ref{fig:predictions} where it can be seen that the error is (linearly) increasing over time.
In conclusion, the computed estimates of the parameters (the learned model) and the initial value yield reliable predictions of the dynamical behaviour of the mass on car system, and thus allow for the application of Funnel-MPC. Future work will be devoted to determine the interplay of learning and Funnel-MPC. To this end, the model will be updated in a receding-horizon fashion to gradually improve the prediction accuracy while keeping the funnel controller as a safeguard.

\section{Conclusion}

In the present paper we proposed the novel method of Funnel-MPC, which ensures initial and recursive feasibility of MPC 
by exploiting concepts from funnel control. To this end, we proposed a new stage cost formulation, which seems to be highly suitable to properly address output-constrained problems with MPC. Moreover, we indicated a further combination of funnel control and MPC by first learning a dynamic model for utilization in MPC, which should lead to improvement of performance and relaxation of the requirements on the sampling rate for a ZOH-based implementation.


\begin{thebibliography}{17}
\providecommand{\natexlab}[1]{#1}
\providecommand{\url}[1]{\texttt{#1}}
\providecommand{\urlprefix}{URL }
\expandafter\ifx\csname urlstyle\endcsname\relax
  \providecommand{\doi}[1]{doi:\discretionary{}{}{}#1}\else
  \providecommand{\doi}{doi:\discretionary{}{}{}\begingroup
  \urlstyle{rm}\Url}\fi

\bibitem[{Berger et~al.(2018)Berger, L{\^e}, and Reis}]{BergLe18}
Berger, T., L{\^e}, H.H., and Reis, T. (2018).
\newblock Funnel control for nonlinear systems with known strict relative
  degree.
\newblock \emph{Automatica}, 87, 345--357.

\bibitem[{Berger et~al.(2019)Berger, Otto, Reis, and Seifried}]{BergOtto19}
Berger, T., Otto, S., Reis, T., and Seifried, R. (2019).
\newblock Combined open-loop and funnel control for underactuated multibody
  systems.
\newblock \emph{Nonlinear Dynamics}, 95(3), 1977--1998.

\bibitem[{Berger and Rauert(2018)}]{BergRaue18}
Berger, T. and Rauert, A.L. (2018).
\newblock A universal model-free and safe adaptive cruise control mechanism.
\newblock In \emph{Proceedings of the MTNS 2018}, 925--932. Hong Kong.

\bibitem[{Berger and Reis(2014)}]{BergReis14a}
Berger, T. and Reis, T. (2014).
\newblock Zero dynamics and funnel control for linear electrical circuits.
\newblock \emph{J. Franklin Inst.}, 351(11), 5099--5132.

\bibitem[{Boccia et~al.(2014)Boccia, Gr{\"u}ne, and Worthmann}]{BoccGrun14}
Boccia, A., Gr{\"u}ne, L., and Worthmann, K. (2014).
\newblock Stability and feasibility of state constrained {MPC} without
  stabilizing terminal constraints.
\newblock \emph{Systems \& control letters}, 72, 14--21.

\bibitem[{Gr{\"u}ne and Worthmann(2008)}]{GrunWort08b}
Gr{\"u}ne, L. and Worthmann, K. (2008).
\newblock Sampled-data redesign for nonlinear multi-input systems.
\newblock In \emph{Geometric Control And Nonsmooth Analysis: In Honor of the
  73rd Birthday of H Hermes and of the 71st Birthday of RT Rockafellar},
  206--227. World Scientific.

\bibitem[{Gr{\"u}ne and Pannek(2017)}]{GruePann17}
Gr{\"u}ne, L. and Pannek, J. (2017).
\newblock \emph{Nonlinear Model Predictive Control}, 45--69.
\newblock Springer International Publishing.

\bibitem[{Gr{\"u}ne et~al.(2008)Gr{\"u}ne, Worthmann, and
  Ne{\v{s}}i{\'c}}]{GrunWort08}
Gr{\"u}ne, L., Worthmann, K., and Ne{\v{s}}i{\'c}, D. (2008).
\newblock Continuous-time controller redesign for digital implementation: a
  trajectory based approach.
\newblock \emph{Automatica}, 44(1), 225--232.

\bibitem[{Hackl(2017)}]{Hack17}
Hackl, C.M. (2017).
\newblock \emph{Non-identifier Based Adaptive Control in Mechatronics--Theory
  and Application}, volume 466 of \emph{Lecture Notes in Control and
  Information Sciences}.
\newblock Springer-Verlag, Cham, Switzerland.

\bibitem[{Ilchmann and Ryan(2008)}]{IlchRyan08}
Ilchmann, A. and Ryan, E.P. (2008).
\newblock High-gain control without identification: a survey.
\newblock \emph{GAMM Mitt.}, 31(1), 115--125.

\bibitem[{Ilchmann et~al.(2002)Ilchmann, Ryan, and Sangwin}]{IlchRyan02b}
Ilchmann, A., Ryan, E.P., and Sangwin, C.J. (2002).
\newblock Tracking with prescribed transient behaviour.
\newblock \emph{ESAIM: Control, Optimisation and Calculus of Variations}, 7,
  471--493.

\bibitem[{Isidori(1995)}]{Isid95}
Isidori, A. (1995).
\newblock \emph{Nonlinear Control Systems}.
\newblock Communications and Control Engineering Series. Springer-Verlag,
  Berlin, 3rd edition.

\bibitem[{Monaco et~al.(2010)Monaco, Normand-Cyrot, and Tiefensee}]{MonaNorm10}
Monaco, S., Normand-Cyrot, D., and Tiefensee, F. (2010).
\newblock Sampled-data redesign of stabilizing feedback.
\newblock In \emph{Proceedings of the 2010 American Control Conference},
  1805--1810. IEEE.

\bibitem[{Rawlings et~al.(2018)Rawlings, Mayne, and Diehl}]{RawlMayn18}
Rawlings, J.B., Mayne, D.Q., and Diehl, M.M. (2018).
\newblock \emph{{Model Predictive Control: Theory, Computation, and Design}}.
\newblock Nob Hill Publishing, second edition.

\bibitem[{Seifried and Blajer(2013)}]{SeifBlaj13}
Seifried, R. and Blajer, W. (2013).
\newblock Analysis of servo-constraint problems for underactuated multibody
  systems.
\newblock \emph{Mech. Sci.}, 4, 113--129.

\bibitem[{Worthmann et~al.(2014)Worthmann, Reble, Gr{\"u}ne, and
  Allg{\"o}wer}]{WortRebl14}
Worthmann, K., Reble, M., Gr{\"u}ne, L., and Allg{\"o}wer, F. (2014).
\newblock The role of sampling for stability and performance in unconstrained
  nonlinear model predictive control.
\newblock \emph{SIAM Journal on Control and Optimization}, 52(1), 581--605.

\bibitem[{Worthmann et~al.(2015)Worthmann, Reble, Gr{\"u}ne, and
  Allg{\"o}wer}]{WortRebl15}
Worthmann, K., Reble, M., Gr{\"u}ne, L., and Allg{\"o}wer, F. (2015).
\newblock Unconstrained nonlinear mpc: Performance estimates for sampled-data
  systems with zero order hold.
\newblock In \emph{2015 54th IEEE Conference on Decision and Control (CDC)},
  4971--4976.

\end{thebibliography}
\end{document}